\documentclass{article}
\usepackage[utf8]{inputenc}
\usepackage[colorinlistoftodos]{todonotes}
\usepackage{subfigure}

\usepackage{amsmath, amssymb, amsthm, graphicx}

\usepackage{etoolbox}

\title{A Reduced Upper Bound for an Edge-coloring Problem from Relation Algebra }

\author{Jeremy F. Alm\\Department of Mathematics\\Illinois College\\Jacksonville, IL 62650\\
\texttt{alm.academic@gmail.com} 
\and David A. Andrews\\Department of Mathematics\\University of Dallas\\Irving, TX 75062\\
\texttt{andrews@udallas.edu}}

\allowdisplaybreaks

\date{\today}

\newtheorem{problem}{Problem}
\newtheorem{theorem}[problem]{Theorem}

\def\M{\mathcal{M}}

\makeatletter
\newcommand{\ct@b}{,\,\discretionary{}{}{}}
\newcounter{listitem@number}
\newcommand{\setinline}[1]{\setcounter{listitem@number}{1}$\{%
  \forcsvlist{\dropin@item}{#1}
  \}$\relax%
}
\newcommand{\dropin@item}[1]{\ifnumequal{\value{listitem@number}}{1}{#1}{$\ct@b$#1}\stepcounter{listitem@number}}
\makeatother

\begin{document}

\maketitle

\begin{abstract}
We construct an edge-coloring of $K_{N}$ (for $N = 3432$) in colors red, dark blue, and light blue, such that there are no monochromatic blue triangles and such that the coloring satisfies a certain strong universal-existential property.  The edge-coloring of $K_{N}$ depends on a cyclic coloring of $K_{17}$ whose two color classes are $K_{4}$-, $K_{4,3}$-, and $K_{5,2}$-free.  This construction yields the smallest known representation of the relation algebra $32_{65}$, reducing the upper bound from 8192 to 3432.
\end{abstract}

\section{Introduction}

In this paper, we consider an edge-coloring problem for complete graphs. Let $K_{N}$ denote the complete graph on $N$ vertices with vertex set $V = V(K_{N})$ and edge set $E = E(K_{N})$. For $x,y \in V$, let $xy$ denote the edge between the two vertices $x$ and $y$. Let $L$ be any finite set and $\M \subseteq L^{3}$. Let $c : E \rightarrow L$.  

For $x,y,z \in V$, let $c(xyz)$ denote the ordered triple $\left(c(xy), c(yz), c(xz)\right)$. We say that $c$ is \emph{good with respect to} $\M$ if the following conditions obtain:
\begin{enumerate}
\item $\forall x,y \in V$ and $\forall (c(xy),j,k) \in \M$, $\exists z \in V$ such that $c(xyz) = (c(xy),j,k)$;\label{condition1}
\item $\forall x,y,z \in V$, $c(xyz) \in \M$; and\label{condition2}
\item $\forall x \in V \ \forall \ell\in L \ \exists \, y\in V$  such that $ c(xy)=\ell $.\label{condition3}
\end{enumerate}

If $K = K_{N}$ has a coloring $c$ which is good with respect to $\mathcal{M}$, then we say that $K$ \emph{realizes} $\M$ (or that $\M$ is \emph{realizable}).

Condition (\ref{condition2}) says that the only triangles allowed are those from $\M$. Condition (\ref{condition1}) says that any triangle that is \emph{allowed} is \emph{mandatory}; in the words of Roger Maddux, ``Anything that can happen, must happen." 

These conditions may seem stringent, but they arise naturally in several situations.  For instance, let  $L=\{r_1, \ldots, r_\ell \}$, and let $$\mathcal{M}_\ell  =\left\{(r_i,r_j,r_k) :|\{i,j,k\}|\in \{ 1,3\} \right\}.$$  Lyndon proved in \cite{lyndon1961} that $\mathcal{M}_\ell$ is realizable on some complete graph if and only if there exists a affine plane of order $\ell -1$, for $\ell >2$. In particular, if $K_N$ realizes $\M_\ell$, then $K_N$ is tiled by edge-disjoint monochromatic $K_{\ell-1}$'s and viewing the vertices of $K_N$ as points and the sets of vertices in the monochromatic $K_{\ell-1}$'s as lines yields an affine plane of order $\ell -1$.

The main motivation for (i)-(iii) is as follows. If we take $R_\alpha=\{(x,y): \ c(xy)= \alpha \}$, and let $ \circ $ stand for ordinary composition of binary relations, i.e., 
\[
  R_\alpha  \circ R_\beta:=\{(x,z): \exists y \ (x,y)\in R_\alpha, (y,z) \in R_\beta\},
\]
then conditions (i) and (ii) imply 
\[
  (R_\alpha  \circ R_\beta) \cap R_\gamma \neq \emptyset \Longrightarrow \ R_\gamma \subseteq R_\alpha \circ R_\beta.
\]
Thus the relations $R_\alpha$ are ``atoms" in the sense that composing them cannot split them into smaller relations.

The problem under consideration is as follows:
\begin{problem}
Find the smallest $N$ such that $\M_n$ is realizable on $K_N$, where $L=\{ r,b_0,...,b_{n-1}\}$ and 
\[
  \M_{n}  = \{r,b_{0},\ldots,b_{n-1}\}^{3}\setminus \{b_{0},\ldots,b_{n-1}\}^{3}.
\]

\end{problem}
 
This problem arises from considerations in algebraic logic, explained in \cite{AMM} and \cite{DH}. Any realization of $\M_n$ gives a representation of a relation algebra; in particular, $\M_2$ gives a representation of the algebra numbered $32_{65}$ in \cite{Madd}.

More generally, given $L$ and $\mathcal{M} \subseteq L^3$, one may ask for which $N$ it is the case that $\mathcal{M}$ is realizable on $K_N$.  In \cite{MadJipTuz}, Maddux, Jipsen and Tuza show that for $\mathcal{M} =L^3$, $K_N$ realizes $\mathcal{M}$ for arbitrarily large finite $N$. In \cite{AndMadd94}, Andr\'eka and Maddux consider all $\M$ when $|L|=2$, and in each case determine the minimum $N$ for which $\M$ is realizable on $K_N$.

It is not difficult to show that for all $n$, $\M_n$ is realizable on $K_\omega$, the complete graph on countably many vertices. It was first shown in \cite{AMM} that $\M_n$ is realizable on $K_N$ for $N<\omega$   using probabilistic methods.  The proof yields a realization of $\M_2$ on $K_N$ for $N$ approximately 7 trillion.

Let 
\[
  f(n) = \min \{ N : \M_n \textrm{ is realizable on } K_N \}.
\]
From \cite{AMM} we have that $f(2) \lessapprox 7 \times 10^9$.  Recently, Dodd and Hirsch \cite{DH}\footnote{The title of \cite{DH} contains a mistake; they improve the \emph{upper} bound.} modified the proof from \cite{AMM} using the Lov\'{a}sz Local Lemma to get $f(2) \lessapprox 37 \times 10^6$.  In June 2014, Dave Sexton and the first author \cite{AS} used direct powers of $\mathbb{Z}_2$, along with a combination of direct construction and computer-assisted randomization, to construct graphs yielding the following bounds:

\begin{align*}
  f(2) &\leq 2^{13} = 8192\\
  f(3) &\leq 2^{16} \\
  f(4) &\leq 2^{19} \\
\end{align*}

In each of \cite{AMM,AS,DH}, the general strategy was as follows: construct large realizations of $\M_1$, and then split the blue edges randomly into two shades of blue, light and dark, using the probabilistic method (in \cite{AMM,DH}) or by  computer-assisted randomization (in \cite{AS}). In the present paper, we improve the upper bound on $f(2)$:

\begin{theorem}\label{thm}
  Let $L=\{ r,b_0,...,b_{n-1}\}$, 
  \[
    \M_{n}  = \{r,b_{0},\ldots,b_{n-1}\}^{3}\setminus \{b_{0},\ldots,b_{n-1}\}^{3},
  \]
  and 
  \[
    f(n) = \min \{ N : \M_n \textrm{ is realizable on } K_N \}.
  \]
  Then $f(2) \leq {14 \choose 7} = 3432$.
\end{theorem}

We achieve better results because our method of splitting the blue edges is entirely deterministic, and is ``uniform" in a certain sense.

\section{Proof of Theorem \ref{thm}}\label{sec:main}

We consider the complete graph $G$ on the size-$7$ subsets of the set of
\emph{points} $[14] = \{1,2,\dots,14\}$. That is, $V = V(G) = \binom{[14]}{7}$. We
provide an edge-coloring that is good with respect to $\M_{2}$ based the size
of intersection of the sets at each vertex. Suppose $X$ and $Y$ are distinct
elements of $\binom{[14]}{7}$.  Define $c$ by
\begin{itemize}
\item if $|X \cap Y| = 0$ then $c(XY) = b_{0}$;
\item if $|X \cap Y| = 1$ then $c(XY) = b_{1}$;
\item if $|X \cap Y| = 2$ then $c(XY) = b_{i}$ for some $i = 0,1$ as described
  below;
\item if $|X \cap Y| \geq 3$ then $c(XY) = r$.
\end{itemize}
A useful notion is that of a \emph{$j$-edge}: if $|X \cap Y| = j$ then we call
$XY$ a $j$-edge. Then the above coloring can be summarized as: all 0-edges are
$b_{0}$, all 1-edges are $b_{1}$, all 3-and-higher-edges are $r$, and all
2-edges are either $b_{0}$ or $b_{1}$ based on a rule to be determined.

Immediately we can show that this coloring has no blue triangles: if $X_{1}$,
$X_{2}$, and $X_{3}$ form an all-blue triangle, then for $i \neq j$, $|X_{i}
\cap X_{j}| \leq 2$ so
\begin{equation*}
  |X_{1} \cup X_{2} \cup X_{3}|
  = \sum_{i} |X_{i}|  - \sum_{i\neq j} |X_{i} \cap X_{j}| + |X_{1} \cap X_{2} \cap X_{3}|  \geq 15,
\end{equation*}
a contradiction since there are only $14$ points.

To determine the color of the 2-edges we will use a \emph{splitting graph}
$G_{s}$. This is a complete graph on the points $[14]$ with all edges colored
either $b_{0}$ or $b_{1}$. We use this coloring to induce a coloring on the
2-edges of $E$: if $X \cap Y = \{i, j\}$ with $i,j$ distinct (so that $XY$ is a
2-edge) then $c(XY)$ will be the color of the edge $ij$ in $G_{s}$.

% \todo[color=green,inline]{The following paragraph needs work!}
% It is important to note that while the coloring of $G_{s}$ is used to color
% some of the edges in $G$, it is a different coloring entirely.  In fact, we
% will see in the following that the edge-coloring of $G_{s}$ needs to satisfy
% three subgraph color properties: no $K_{4}$ subgraph, $K_{3,3}$ subgraph, nor
% $K_{5,2}$ subgraph of $G_{s}$ can be monochromatic.

It is important to note that while the coloring of $G_{s}$ is used to color
some of the edges in $G$, specifically the 2-edges, the graphs $G$ and $G_{s}$
are entirely different graphs. In the argument below we will need a coloring of
the edges of $G_{s}$ that contains no monochromatic $K_{4}$, no monochromatic
$K_{4,3}$, and no monochromatic $K_{5,2}$. Graphs with all three of these
properties do exist. For example, we consider the cyclic coloring of
$K_{17}$ in which one of the color classes (which we will color $b_{0}$)
is given by $\mathcal{C}_{0} = \{1,2,4,8,9,13,15,16\}$ and the other (which we
will color $b_{1}$) is given by $\mathcal{C}_{1} = \{3,5,6,7,10,11,12,14\}$.
The authors verified that this coloring has no monochromatic $K_{4}$, $K_{4,3}$,
or $K_{5,2}$ subgraphs, and so for $G_{s}$ we could use the subgraph obtained
from this $K_{17}$ by deleting vertices 15, 16, and 17. (We verified this using two different programs written in two different languages, \texttt{Java} and \texttt{Python}.
The authors are amazed that $R(4,4) = R(K_{3,3},K_{3,3}) = R(K_{5,2},K_{5,2}) = 18$ and the lower bound for all three Ramsey numbers can be established by the same graph).

To see how the coloring of $G_{s}$ determines the coloring of edges in $G$, consider the vertices
\begin{align*}
  X &=\{1, 2, 3, 4, 5, 6, 7\},\\
  Y &=\{2, 3, 4, 5, 6, 7, 8\}, \textrm{ and}\\
  Z &=\{7, 8, 9, 10, 11, 12, 13\}.
\end{align*}
and the coloring of the edges of the triangle $XYZ$. Since $|X\cap Y| = 6 \geq 3$, $c(XY)=r$. Since
$|X\cap Z| = 1$, $c(XZ)=b_1$. Since $|Y\cap Z| = 2$, we must look at $G_s$. We
have $Y\cap Z=\{7,8\}$, so we look at the color of the edge between
the vertices labeled 7 and 8 in the example $G_{s}$ given above. Since $8-7=1\in\mathcal{C}_0$, $c(YZ)=b_0$.

We introduce a bit of terminology that will prove to be convenient. Given an edge $XY$ colored $c_1$, if $(c_1,c_2,c_3)\in \M_2$, we will say that $XY$ has the \emph{need} $(c_1,c_2,c_3)$ (as required by (i)). If $Z$ is such that $c(XZ)=c_2$ and $c(YZ)=c_3$, we will call $Z$ a \emph{witness} to the need $(c_1,c_2,c_3)$. Note that, because there can be no all-blue triangles (by construction), (ii) is automatically satisfied; thus in order to show that $K_N$ is a realization of $\M_2$, we must show that every edge has its needs met, i.e., condition (i) is satisfied.  (Note that, in the case of $\M_2$, (i)$\Rightarrow$(iii),  so we need not concern ourselves with (iii).)

To show that each edge in $G$ has all its needs met, we will consider without
loss of generality only edges between the vertex $X =$ \setinline{1,2,3,4,5,6,7}
and a vertex $Y$ that overlaps $X$ in its last $j$ elements. For example, for
$j=3$ we will consider
$X =$ \setinline{1,2,3,4,5,6,7} and
$Y =$ \setinline{5,6,7,8,9,10,11}.
Furthermore, in order not to lose any generality in the proof, we will not
assume anything about the particular coloring of $G_{s}$ other than the
monochromatic-subgraph-free properties given above.

For $j=0,1,$ and $2$, the edge $XY$ is colored either $b_{0}$ or $b_{1}$, and
has five needs. If $c(XY)=b_{i}$ these needs are $(b_{i},r,b_{0})$,
$(b_{i},r,b_{1})$, $(b_{i},b_{0},r)$ and $(b_{i},b_{1},r)$, and $(b_{i},r,r)$.
For $j=3,4,5,$ and $6$, the edge $XY$ is colored $r$ and has 9 needs, of four
different categories:
\begin{itemize}
\item the \emph{homogeneous blue} needs: $(r,b_{0},b_{0})$ and $(r, b_{1}, b_{1})$;
\item the \emph{heterogeneous blue} needs: $(r,b_{0},b_{1})$ and $(r,b_{1},b_{0})$;
\item the \emph{red-blue} needs: $(r,r,b_{i})$ and $(r,b_{i},r)$ for $i \in
  \{0,1\}$; and
\item the \emph{all red} need: $(r,r,r)$.
\end{itemize}

\subsection{Notation}

Throughout the remainder of Section \ref{sec:main}, we will exhibit, for every edge $XY$ and for each need $(c_1,c_2,c_3)$, a witness $Z$ satisfying that need.  It will be convenient to introduce the following piece of notation:
\[
 \{1,2,3,4,5,[8,9,10,11]_{2}\}
\]
will stand for the collection of six sets
\begin{gather*}
  \{1,2,3,4,5,8,9\},\ \{1,2,3,4,5,8,10\},\ \{1,2,3,4,5,8,11\}, \\
  \{1,2,3,4,5,9,10\},\ \{1,2,3,4,5,9,11\},\ \{1,2,3,4,5,10,11\},
\end{gather*}
where each set contains 1, 2, 3, 4, 5, and exactly two from 8, 9, 10, and 11. We will demonstrate the usefulness of this notation in the next section.

% Suppose, for example, that $X=\{1,2,3,4,5,6,7\}$ and $Y=\{8,9,10,11,12,13,14\}$, so that $c(XY)=b_0$. Choosing a vertex $Z$ from $\{1,2,3,4,5,[8,9,10,11]_{2}\}$ will yield $|X\cap Z|=5$, so $c(XZ)=r$, and $|Y\cap Z|=2$, so $c(YZ)=b_i$, where $b_i$ is the color in $G_s$ of the edge between the two points in $Y\cap Z$.  Since, by assumption, $G_{s}$ contains no monochromatic $K_4$, we may choose $Z_0$ and $Z_1$ from $\{1,2,3,4,5,[8,9,10,11]_{2}\}$ such that $c(YZ_0)=b_0$ and $c(YZ_1)=b_1$. So in this case, we would say that $\{1,2,3,4,5,[8,9,10,11]_{2}\}$ witnesses the needs $(b_0, r, b_0)$ and $(b_0, r, b_1)$ for $XY$.

\subsection{$0$-edge case}

We work out the details carefully in this case as the other cases use similar
approaches.

For this case, we are considering $X =$ \setinline{1,2,3,4,5,6,7} and
$Y =$ \setinline{8,9,10,11,12,13,14}.
By our coloring $c(XY) =b_{0}$ and so we must show
there are vertices that witness the five needs: $(b_{0},r,b_{0})$,
$(b_{0},r,b_{1})$, $(b_{0},b_{0},r)$ and $(b_{0},b_{1},r)$, and $(b_{0},r,r)$.

To show that the first need is satisfied,
% it suffices to find a vertex $Z_{1} \in V$ different from $X$ and $Y$ such
% that $c(XZ_{1}) = r$ and $c(YZ_{1}) = b_{0}$.  The first requirement is easy
% to satisfy, as we must simply make sure that $|X \cap Z_{1}| \geq 3$. There
% are two ways to get $c(YZ_{1}) = b_{0}$: we could make $Y$ and $Z_{1}$
% disjoint, or we could make $Y \cap Z_{1} = \{i,j\}$ such that color of the
% $ij$ edge in $G_{s}$ is $b_{0}$.  In this case it is impossible to make $Y$
% and $Z_{1}$ disjoint, since $X$ is the only such set, so we must find an
% appropriate $i,j$ combination to put into $Z_{1}$.
we consider selecting $Z$ from the six-set collection
\setinline{1,2,3,4,5,[8,9,10,11]_{2}}. All choices from this collection overlap $X$
in 5 points --- giving $c(XZ) = r$ --- and overlap $Y$ in 2 points so that the
edge $YZ$ is colored the same as the corresponding edge in $G_{s}$. The
subgraph of $G_{s}$ induced by the points $\{8,9,10,11\}$ form a $K_{4}$
subgraph, and so is not monochromatic. Thus there must be an edge $ij$ in
this subgraph that is colored $b_{0}$. Thus we can use $Z=$
\setinline{1,2,3,4,5,i,j} to
witness this first need. Similarly, there must be an edge $k\ell$ in this subgraph
that has color $b_{1}$, and so $Z' =$ \setinline{1,2,3,4,5,k,ell}
that satisfies the second need: $c(XZ') = r$ and $c(YZ') = b_{1}$.

A similar construction can be used to find vertices that witness the
$(b_{0},b_{i},r)$ needs.  The need $(b_{0},r,r)$ is also witnessed, and so we
summarize:
\begin{itemize}
\item $(b_{0},r,b_{0})$ and $(b_{0},r,b_{1})$ have witnesses from
  \setinline{1,2,3,4,5,[8,9,10,11]_{2}}; 
\item $(b_{0},b_{0},r)$ and $(b_{0},b_{1},r)$ have witnesses from
  \setinline{[1,2,3,4]_{2},8,9,10,11,12}; and
\item $(b_{0},r,r)$ is satisfied by the witness \setinline{5,6,7,8,9,10,14}.
\end{itemize}
Thus all needs of all $0$-edges are met in our graph coloring.

\subsection{$1$-edge case}

We consider $X =$ \setinline{1,2,3,4,5,6,7} and $Y =$ \setinline{7,8,9,10,11,12,13},
so that $c(XY) = b_{1}$. In this case, we can provide
witnesses for all 5 needs by direct construction:
\begin{itemize}
\item $(b_{1},r,b_{0})$ is witnessed by \setinline{1,2,3,4,5,6,14} (since $|X\cap Z|=6$ and $|Y\cap Z|=0$);
\item $(b_{1},r,b_{1})$ is witnessed by \setinline{1,2,3,4,5,6,8} (since $|X\cap Z|=6$ and $|Y\cap Z|=1$);
\item $(b_{1},b_{0},r)$ is witnessed by \setinline{8,9,10,11,12,13,14} (since $|X\cap Z|=0$ and $|Y\cap Z|=6$);
\item $(b_{1},b_{1},r)$ is witnessed by \setinline{6,8,9,10,11,12,13} (since $|X\cap Z|=1$ and $|Y\cap Z|=6$); and
\item $(b_{1},r,r)$ is again witnessed by \setinline{5,6,7,8,9,10,14} (since $|X\cap Z|=3$ and $|Y\cap Z|=4$).
\end{itemize}

\subsection{$2$-edge case}

For this case, we consider $X =$ \setinline{1,2,3,4,5,6,7} and $Y =
\{6,7,8,9,10,11,12\}$. We know $c(XY) = b_{i}$ for either $i=0$ or $1$, but
which one is not relevant, as the form of the needs do not differ, and can
be satisfied again by direct construction:
\begin{itemize}
\item $(b_{i},r,b_{0})$ is witnessed by \setinline{1,2,3,4,5,13,14};
\item $(b_{i},r,b_{1})$ is witnessed by \setinline{1,2,3,4,5,8,14};
\item $(b_{i},b_{0},r)$ is witnessed by \setinline{8,9,10,11,12,13,14};
\item $(b_{i},b_{1},r)$ is witnessed by \setinline{7,8,9,10,11,12,13}; and
\item $(b_{i},r,r)$ is also witnessed by \setinline{5,6,7,8,9,10,14}.
\end{itemize}

\subsection{$3$-edge case}

Since $X =$ \setinline{1,2,3,4,5,6,7} and $Y =$ \setinline{5,6,7,8,9,10,11} we have $c(XY) =
r$.  The red-blue needs are satisfied in a similar way as they were in the
$0$-edge case: the $(r,r,b_{i})$ needs are witnessed by some vertices selected
from \setinline{1,2,3,[8,9,10,11]_{2},13,14}, and the $(r,b_{i},r)$ are witnessed by
vertices in \setinline{[1,2,3,4]_{2},9,10,11,13,14}.

The homogeneous and heterogeneous blue needs all have witnesses in the 36 sets in
\setinline{[1,2,3,4]_{2}, [8,9,10,11]_{2}, 12, 13, 14}: the subgraph of $G_{s}$
induced by $\{1,2,3,4\}$ has $b_{0}$ and $b_{1}$ edges, as does the subgraph
induced by $\{8,9,10,11\}$. Furthermore, edges from both of these subgraphs can
be selected independently since the first  involves points from $X$ only and the
second  points from $Y$ only. Thus, we can find witnesses for all of these
needs.

In summary, we have:
\begin{itemize}
\item $(r,b_{i},b_{j})$ for $i,j \in \{0,1\}$) have witnesses in \setinline{[1,2,3,4]_{2}, [8,9,10,11]_{2}, 12, 13, 14};
\item $(r,r,b_{i})$ have witnesses in $\{1,2,3,[8,9,10,11]_{2},12,13\}$;
\item $(r,b_{k},r)$ have witnesses in $\{[1,2,3,4]_{2},8,9,10,12,13\}$; and
\item $(r,r,r)$ has witness $\{5,6,7,8,9,10,14\}$.
\end{itemize}

\subsection{$4$-edge case}

For this case, $X =$ \setinline{1,2,3,4,5,6,7} and $Y = \{4,5,6,7,8,9,10\}$, and $c(XY)
= r$.

We begin with the easiest needs, red-blue and all red:
\begin{itemize}
\item The need $(r,r,b_{0})$ is witnessed by $\{1,2,3,11,12,13,14\}$;
\item The need $(r,r,b_{1})$ is witnessed by $\{1,2,3,10,11,12,13\}$;
\item The need $(r,b_{0},r)$ is witnessed by $\{8,9,10,11,12,13,14\}$;
\item The need $(r,b_{1},r)$ is witnessed by $\{1,9,10,11,12,13,14\}$; and
\item The need $(r,r,r)$ is witnessed by $\{5,6,7,8,9,10,14\}$.
\end{itemize}

To satisfy the heterogeneous and homogeneous blue needs, we need to argue by
cases, which will be fairly complex, but use similar approaches. Thus, we will
rehearse some parts of the argument first.

In these cases, we must rely on $G_{s}$ to generate the correct (blue) coloring
in $b_{i}$ since each vertex $Z$ in $G$ that forms a $0$- or $1$-edge with $X$
must form a $3$-or-higher edge with $Y$ (and similarly for vertices that form a
$0$- or $1$-edge with $Y$). Thus we look at a subgraph of $G_{s}$ with vertex set $X
\cup Y$ arranged in three groups --- a ``spine'' of the
points in $S = X \cap Y$ and two ``wings'' of the points in $W_{X} = X\setminus
Y$ and $W_{Y}=Y \setminus X$ --- and with edges from the union of:
\begin{itemize}
\item the complete graph induced by $W_{X}$;
\item the complete graph induced by $W_{Y}$;
\item the complete bipartite graph induced by the part sets $W_{X}$ and $S$; and
\item the complete bipartite graph induced by the part sets $W_{Y}$ and $S$.
\end{itemize}

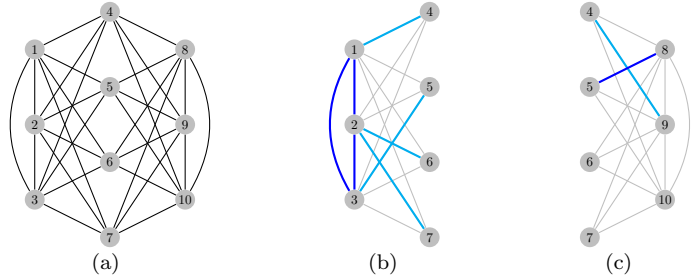
\begin{figure}[ht!]
  \centering
  
  \subfigure[]{%
    \label{fig:s-a-w-overview}%
    \begin{tikzpicture}
      [auto,node distance=1cm,point/.style={circle,minimum
        width=15pt,fill=black!25,inner sep=0pt,text=black,scale=0.5}]
        
      \draw[use as bounding box] (-2,0) rectangle (2,3) [white];
      \node[point] (4) at (0,3) {$4$};
      \node[point] (5) at (0,2) {$5$};
      \node[point] (6) at (0,1) {$6$};
      \node[point] (7) at (0,0) {$7$};

      \node[point] (1) at (-1,2.5) {$1$};
      \node[point] (2) at (-1,1.5) {$2$};
      \node[point] (3) at (-1,0.5) {$3$};

      \node[point] (8) at (1,2.5) {$8$};
      \node[point] (9) at (1,1.5) {$9$};
      \node[point] (10) at (1,0.5) {$10$};

      \path
      \foreach \source in {4,...,7}{
        \foreach \dest in {1,2,3,8,9,10}{
          (\source) edge (\dest)
        }
      };

      \path (1) edge (2)
            (2) edge (3)
            (1) edge [bend right] (3);

      \path (8) edge (9)
            (9) edge (10)
            (8) edge [bend left] (10);
    \end{tikzpicture}}
  \subfigure[]{%
    \label{fig:s-a-w-k4need}
    \begin{tikzpicture}
      [auto,node distance=1cm,point/.style={circle,minimum
        width=15pt,fill=black!25,inner sep=0pt,text=black,scale=0.5}]
      \draw[use as bounding box] (-2,0) rectangle (1,3) [white];
      \node[point] (4) at (0,3) {$4$};
      \node[point] (5) at (0,2) {$5$};
      \node[point] (6) at (0,1) {$6$};
      \node[point] (7) at (0,0) {$7$};

      \node[point] (1) at (-1,2.5) {$1$};
      \node[point] (2) at (-1,1.5) {$2$};
      \node[point] (3) at (-1,0.5) {$3$};

      \path
      \foreach \source in {4,...,7}{
        \foreach \dest in {1,2,3}{
          (\source) edge [thin,black!25] (\dest)
        }
      };
      
      \path (1) edge [thick,blue] (2)
            (2) edge [thick,blue] (3)
            (1) edge [thick,blue,bend right] (3);

      \path (4) edge [thick,cyan] (1)
            (5) edge [thick,cyan] (3)
            (6) edge [thick,cyan] (2)
            (7) edge [thick,cyan] (2);
            
    \end{tikzpicture}}%
  \subfigure[]{%
    \label{fig:s-a-w-k43need}
    \begin{tikzpicture}
      [auto,node distance=1cm,point/.style={circle,minimum
        width=15pt,fill=black!25,inner sep=0pt,text=black,scale=0.5}]
      \draw[use as bounding box] (-1,0) rectangle (2,3) [white];
      \node[point] (4) at (0,3) {$4$};
      \node[point] (5) at (0,2) {$5$};
      \node[point] (6) at (0,1) {$6$};
      \node[point] (7) at (0,0) {$7$};

      \node[point] (8) at (1,2.5) {$8$};
      \node[point] (9) at (1,1.5) {$9$};
      \node[point] (10) at (1,0.5) {$10$};

      \path
      \foreach \source in {4,...,7}{
        \foreach \dest in {8,9,10}{
          (\source) edge [thin,black!25] (\dest)
        }
      };
      
      \path (8) edge [thin,black!25] (9)
            (9) edge [thin,black!25] (10)
            (8) edge [thin,black!25,bend left] (10);

      \path (4) edge [thick,cyan] (9)
            (5) edge [thick,blue] (8);    
    \end{tikzpicture}}
  \caption{(a) The ``spine and wing'' diagram used in the $4$-edge case (a
            subgraph of $G_{S}$).
           (b) Each node of the spine has a edge of the opposite color to a
           monochromatic wing.
           (c) There exists at least one edge in each color from the spine to a wing.
           Colorings in these diagrams are for demonstration purposes only. Actual results may vary. Past performance is not an indicator of future results.}
  \label{fig:spine-and-wings}
\end{figure}

To satisfy each need, we have two mutually exclusive options:
\begin{itemize}
\item choose two points $i,j \in W_{X}$ and two points from $m,n \in W_{Y}$
  with edge $ij$ having the desired color of $XZ$ and $mn$ having the desired
  color of $YZ$; or
\item choose one point from each of $W_{X}$, $S$, and $W_{Y}$ ($i$, $j$, and
  $k$ respectively) where $ij$ has the desired color of $XZ$ and $jk$ has the
  desired color of $YZ$ (hence the ``wing'' nomenclature).
\end{itemize}
The first option can be used to satisfy all the needs trivially if the graphs
induced by $W_{X}$ and $W_{Y}$ each contain edges of both colors (Case I
below).  On the other hand, if one of these graphs is monochromatic (Cases II and III),
the first option can  satisfy only some of the needs.  Thus, we will have to rely on the
second option to satisfy the remaining needs, and it is not  obvious
that this can be done. However, we can show that these needs are satisfied
by using the properties of $G_{s}$ in two particular ways.

First, any point $k \in S$ can be combined with $W_{X}$ (or $W_{Y}$) to form a
$K_{4}$ subgraph of $G_{s}$, which we know is not monochromatic. But if $W_{X}$
is monochromatic, we must then have an edge of the other color from $k$
to one of the points in $W_{X}$. This means that, when one of the wings
is monochromatic, every point in $S$ has an edge of the other color to some
point in $W_{X}$ (see Figure \ref{fig:s-a-w-k4need}).

Second, we know that the $K_{3,4}$ graph induced by $L = S$ and $R = W_{Y}$ (or $W_{X}$)
is not monochromatic, and so there must be an edge of each color
from $S$ to $W_{Y}$ (see Figure \ref{fig:s-a-w-k43need}).

\begin{description}
\item[Case I:] the graphs induced by $W_{X}$ and $W_{Y}$ are not
  monochromatic. Then each blue-only need is witnessed by a vertex in
  $\{[1,2,3]_{2},[8,9,10]_{2},11,12,13\}$.
\item[Case II:] the graph induced by $W_{X}$ is monochromatic $b_{0}$.  
  \begin{description}
  \item[Case II.A:] $W_{Y}$ is not monochromatic. In this case there are
    vertices in $\{1,2,[8,9,10]_{2},12,13,14\}$ that witness to the needs
    $(r,b_{0},b_{i})$ for $i=0,1$. 

    We consider the $K_{4,3}$ subgraph induced by $S$ and
    $W_{Y}$: there must be an edge $ij$ ($i \in S$, $j \in W_{Y}$)
    that is $b_{0}$ (as in Figure~\ref{fig:spine-and-wings}(c)). Furthermore,
    in the subgraph induced by $\{1,2,3,i\}$, there must be an edge $ik$
    ($k \in W_{X}$) that is $b_{1}$ (as in Figure~\ref{fig:spine-and-wings}(b)).
    We now have a witness to the $(r,b_{1},b_{0})$ need: $\{k,i,j,11,12,13,14\}$.

    In a similar way, there must be an edge $mn$ ($m \in S$, $n \in W_{Y}$)
    that is $b_{1}$ and an edge $\ell{}m$ in the subgraph induced by
    $\{1,2,3,m\}$ ($\ell \in W_{X}$) that has color $b_{1}$, and so $\{\ell, m,
    n, 11, 12, 13, 14\}$ is a witness for $(r,b_{1},b_{1})$.
    
  \item[Case II.B:] $W_{Y}$ is monochromatic $b_{0}$. We immediately have a
    witness for $(r, b_{0}, b_{0})$: $\{1,2,9,10,11,12,13\}$.

    We consider the bipartite subgraph induced by $S$ and
    $W_{Y}$, which must have an edge $ij$ ($i \in S$, $j \in W_{Y}$)
    that is $b_{1}$; and the subgraph induced by $\{1,2,3,i\}$ must have and
    edge $i\ell$ ($\ell \in W_{X}$) that is $b_{1}$.  Thus $\{\ell, i, j,
    11, 12, 13, 14\}$ is a witness for $(r, b_{1}, b_{1})$.

    The bipartite graph in the previous paragraph also has an edge edge $mn$
    ($m \in S$, $n \in W_{Y}$) that is $b_{0}$; and the subgraph induced by
    $\{1,2,3,m\}$ must have an edge $mk$ ($k \in W_{X}$) that is $b_{1}$. Thus
    $\{k, m, n, 11, 12, 13, 14\}$ is a witness for $(r, b_{1}, b_{0})$.

    Finally, we consider another bipartite subgraph, the one induced by $S$ and
    $W_{X}$. There must be an edge $pq$ ($p \in W_{X}$, $q \in S$) that is
    $b_{0}$. Since the subgraph induced by $\{8,9,10,q\}$ is not monochromatic
    (but $W_{Y}$ is monochromatic $b_{0}$) there must be an edge $qs$ ($s \in
    W_{Y}$) that is $b_{1}$. Thus $\{p, q, s, 11, 12, 13, 14\}$ is a witness
    for $(r, b_{0}, b_{1})$.
    
  \item[Case II.C:] $W_{Y}$ is monochromatic $b_{1}$. In this case, we
    immediately have $\{1,2,9,10,11,12,13,14\}$ as a witness for
    $(r,b_{0},b_{1})$.

    Once again, we consider the bipartite subgraph induced by $S$ and
    $W_{Y}$. There must be an edge $ij$ ($i \in S$, $j \in W_{Y}$) that is
    $b_{1}$. But considering the subgraph induced by $\{1,2,3,i\}$, there must
    be an edge $i\ell$ ($\ell \in W_{X}$) that is $b_{1}$. Thus $\{\ell, i,
    j, 11, 12, 13, 14\}$ is a witness of $(r, b_{1}, b_{1})$.

    We can satisfy another need by looking at the bipartite graph induced by
    $S$ and $W_{X}$: there must be an edge $mn$ ($m \in S$, $n \in W_{X}$) that
    is $b_{0}$.  However, in the subgraph induced by $\{8,9,10,m\}$ there must
    be an edge $mk$ ($k \in W_{Y}$) that is $b_{0}$. And so, $\{k, m, n, 11, 12,
    13, 14\}$ is a witness to $(r, b_{0}, b_{0})$.

    For the final need we take a slightly different approach and first consider
    the subgraph induced by $\{1,2,3,4\}$ (4 was chosen arbitrarily from $S$
    --- any point of $S$ would do). Since this $K_{4}$ subgraph is not
    monochromatic (but $W_{X}$ is monochromatic $b_{0}$) there is an $p \in
    W_{X}$ such that edge $4p$ is $b_{1}$. Similarly, the subgraph induced by
    $\{4,8,9,10\}$ must have an edge $4q$ ($q \in W_{Y}$) that is
    $b_{0}$. Thus, $\{p, q, 4, 11, 12, 13, 14\}$ is a witness of
    $(r,b_{1},b_{0})$.
  \end{description}
\item[Case III:] the graph induced by $W_{X}$ is monochromatic $b_{1}$. This
  proceeds in the same way as the previous case.
\end{description}

\subsection{$5$-edge case}

For this case, $X =$ \setinline{1,2,3,4,5,6,7} and $Y = \{3,4,5,6,7,8,9\}$, and $c(XY)
= r$.

It is straight-forward to see that the red-blue, homogeneous blue, and all-red needs are
satisfied:
\begin{itemize}
\item $(r,r,b_{i})$ have witnesses in $\{1,2,[3,4,5,6]_{2},12,13,14\}$;
\item $(r,b_{k},r)$ have witnesses in $\{[3,4,5,6]_{2},8,9,12,13,14\}$; and
\item $(r,b_{0},b_{0})$ and $(r,b_{1},b_{1})$ are both satisfied by witnesses
  in \setinline{[3,4,5,6]_{2},10,11,12,13,14}.
\item $(r,r,r)$ is witnessed by $\{5,6,7,8,9,10,14\}$.
\end{itemize}

How the heterogeneous blue needs --- $(r,b_{0},b_{1})$ and $(r,b_{1},b_{0})$ ---
are satisfied depends on the coloring of the 1--2 and 8--9 edges in $G_{s}$.
If both edges are the same color, we can construct witnesses for each of the
needs, as follows:
\begin{itemize}
\item if both 1--2 and 8--9 are $b_{0}$, then $\{1,2,8,11,12,13,14\}$
  witnesses $(r,b_{0},b_{1})$ and $\{1,8,9,11,12,13,14\}$ witnesses $(r,b_{1},b_{0})$;
\item if both 1--2 and 8--9 are $b_{1}$, then $\{8,9,10,11,12,13,14\}$
  witnesses $(r,b_{0},b_{1})$ and $\{1,2,10,11,12,13,14\}$ witnesses $(r,b_{1},b_{0})$;.
\end{itemize}

If 1--2 is $b_{0}$ and 8--9 is $b_{1}$ then $\{1,2,8,9,12,13,14\}$ is a
witness of $(r,b_{0},b_{1})$. Now consider the vertices $Z$ in
$\{[3,4,5,6,7]_{1}, [8,9]_{1}, 10, 11, 12, 13, 14\}$; for all of these $c(XZ) = b_{1}$, but
$c(YZ)$ depends on the points selected. However, the bipartite subgraph
induced by the sets $L=\{3,4,5,6,7\}$ and $R=\{8,9\}$ is $K_{5,2}$ and so has
one edge $ij$ that is colored $b_{0}$, and so $\{i,j,10,11,12,13,14\}$
witnesses $(r,b_{1},b_{0})$. A very similar argument works if 1--2 is $b_{1}$
and 8--9 is $b_{0}$.

\subsection{$6$-edge case}

Our final case considers $X =$ \setinline{1,2,3,4,5,6,7} and $Y = \{2,3,4,5,6,7,8\}$.
The heterogeneous blue and all-red needs are satisfied by directly constructed
witnesses:
\begin{itemize}
\item $\{1,9,10,11,12,13,14\}$ is a witness of $(r,b_{1},b_{0})$;
\item $\{8,9,10,11,12,13,14\}$ is a witness of $(r,b_{0},b_{1})$;
\item $\{5,6,7,8,9,10,14\}$ is a witness of $(r,r,r)$.
\end{itemize}
The homogeneous blue needs --- $(r,b_{0},b_{0})$ and $(r,b_{1},b_{1})$ --- must
each have a witness in $\{[2,3,4,5]_{2},10,11,12,13,14\}$.  Finally the needs
$(r,r,b_{i})$ are witnessed by some vertices in
$\{1,[2,3,4,5]_{2},11,12,13,14\}$ and the needs $(r,b_{i},r)$ are witnessed by
vertices in $\{[2,3,4,5]_{2},8,11,12,13,14\}$.

\section{Conclusion}

It is extremely doubtful that 3432 is the correct value of $f (2)$; the authors guess that  $f (2)<1000$.  However, further improvement would seem to require the construction of triangle-free graphs that are both dense and ``uniform" (vertex-transitive, for example).  One natural approach would be to use cyclic  colorings, where the ``blue" color class is given by a maximal sum-free subset of $\mathbb{Z}/N\mathbb{Z}$.  The authors have attempted to construct such a subset, but have failed to find a sum-free subset that has enough ``redundancy", so that once the blue color class is split into light blue and dark blue, all needs are still met.  It might be possible to construct such a sum-free set over a very large modulus, but (after much unfruitful effort) the authors are doubtful of the existence of such a set over a modulus less than 3432.

% \bibliographystyle{amsplain}
% \bibliography{refs}

\providecommand{\bysame}{\leavevmode\hbox to3em{\hrulefill}\thinspace}
\providecommand{\MR}{\relax\ifhmode\unskip\space\fi MR }
\providecommand{\MRhref}[2]{%
  \href{http://www.ams.org/mathscinet-getitem?mr=#1}{#2}
}
\providecommand{\href}[2]{#2}

\end{document}